\documentstyle[12pt, leqno, amstex]{amsart}

\makeatletter

\begin{document}
\bibliographystyle{amsalpha}

\newtheorem{Assumption}{Assumption}[section]
\newtheorem{Theorem}{Theorem}[section]
\newtheorem{Lemma}{Lemma}[section]
\newtheorem{Remark}{Remark}[section]
\newtheorem{Corollary}{Corollary}[section]
\newtheorem{Conjecture}{Conjecture}[section]
\newtheorem{Proposition}{Proposition}[section]
\newtheorem{Example}{Example}[section]
\newtheorem{Definition}{Definition}[section]
\newtheorem{Question}{Question}[section]
\renewcommand{\thesubsection}{\it}

\baselineskip=14pt
\addtocounter{section}{-1}

\title{Symmetries in the system of type $D_4^{(1)}$}

\author{Yusuke Sasano}

\address{Graduate School of Mathematical 
Sciences, The University of Tokyo, Tokyo 153-8914, Japan.}
\email{sasano@ms.u-tokyo.ac.jp}
\keywords{Affine Weyl group, birational symmetries, coupled Painlev\'e systems.}
\thanks{2000 Mathematics Subject Classification Numbers. 34M55, 34M45, 58F05, 32S65.}

\begin{abstract}
In this paper, we propose a 4-parameter family of coupled Painlev\'e III systems in dimension four with affine Weyl group symmetry of type $D_4^{(1)}$. We also propose its symmetric form in which the $D_4^{(1)}$-symmetries become clearly visible.
\end{abstract}
\maketitle

\section{Statement of main results} 
In \cite{Sasa4,Sasa5,Sasa6}, we presented some types of coupled Painlev\'e systems with various affine Weyl group symmetries. In this paper, we present a 4-parameter family of coupled Painlev\'e III systems with affine Weyl group symmetry of type $D_4^{(1)}$. Since these universal B{\"a}cklund transformations have Lie theoretic origin, similarity reduction of a Drinfeld-Sokolov hierarchy admits such a B{\"a}cklund symmetry. The aim of this paper is to introduce the system of type $D_4^{(1)}$. After our discovery of this system, they were studied from the viewpoint of Drinfeld-Sokolov hierarchy by K. Fuji independently (cf. [1]), and he succeeded to obtain our system by similarity reduction of the Drinfeld-Sokolov hierarchy of type $D_4^{(1)}$.

At first, we propose a 4-parameter family of autonomous ordinary differential systems with the invariant divisors $f_i$ as variables:
\begin{equation}
  \left\{
  \begin{aligned}
   \frac{df_0}{dt} &=-(2f_1g_1+\alpha_1)f_0-\alpha_0f_1,\\
   \frac{df_1}{dt} &=-(2f_0g_1+\alpha_0)f_1-\alpha_1f_0,\\
   \frac{df_2}{dt} &=\{(f_0+f_1)g_1+(f_3+f_4)g_2+1\}f_2-2\alpha_2 g_1g_2,\\
   \frac{df_3}{dt} &=-(2f_4g_2+\alpha_4)f_3-\alpha_3f_4,\\
   \frac{df_4}{dt} &=-(2f_3g_2+\alpha_3)f_4-\alpha_4f_3,\\
   \frac{dg_1}{dt} &=(f_0+f_1)g_1^2-\{(f_3+f_4)g_2-\alpha_0-\alpha_1\}g_1+(f_3+f_4)f_2,\\
   \frac{dg_2}{dt} &=(f_3+f_4)g_2^2-\{(f_0+f_1)g_1-\alpha_3-\alpha_4\}g_2+(f_0+f_1)f_2.\\
   \end{aligned}
  \right. 
\end{equation}
Here $f_0,f_1,\dots,f_4$ and $g_1,g_2$ denote unknown complex variables and $\alpha_0,\dots,\alpha_4$ are the parameters satisfying the condition:
$$
\alpha_0+\alpha_1+2\alpha_2+\alpha_3+\alpha_4=1.
$$

\begin{Proposition}
This system has the following invariant divisors\rm{:\rm}
\begin{center}
\begin{tabular}{|c|c|c|} \hline
invariant divisors & parameter's relation \\ \hline
$f_0:=0$ & $\alpha_0=0$  \\ \hline
$f_1:=0$ & $\alpha_1=0$  \\ \hline
$f_2:=0$ & $\alpha_2=0$  \\ \hline
$f_3:=0$ & $\alpha_3=0$  \\ \hline
$f_4:=0$ & $\alpha_4=0$  \\ \hline
\end{tabular}
\end{center}
\end{Proposition}

\begin{Theorem}
This system is invariant under the transformations $s_0,\dots,s_4$\\
defined as follows$:$ with {\it the notation} $(*):=(f_0,f_1,\dots,f_4,g_1,g_2;\alpha_0,\alpha_1,\dots,\alpha_4),$
\begin{align}
\begin{split}
s_0:(*) \rightarrow &(f_0,f_1,f_2+\frac{\alpha_0 g_2}{f_0},f_3,f_4,g_1+\frac{\alpha_0}{f_0},g_2;-\alpha_0,\alpha_1,\alpha_2+\alpha_0,\alpha_3,\alpha_4),\\
s_1:(*) \rightarrow &(f_0,f_1,f_2+\frac{\alpha_1 g_2}{f_1},f_3,f_4,g_1+\frac{\alpha_1}{f_1},g_2;\alpha_0,-\alpha_1,\alpha_2+\alpha_1,\alpha_3,\alpha_4),\\
s_2:(*) \rightarrow &(f_0-\frac{\alpha_2 g_2}{f_2},f_1-\frac{\alpha_2 g_2}{f_2},f_2,f_3-\frac{\alpha_2 g_1}{f_2},f_4-\frac{\alpha_2 g_1}{f_2},g_1,g_2;\\
&\alpha_0+\alpha_2,\alpha_1+\alpha_2,-\alpha_2,\alpha_3+\alpha_2,\alpha_4+\alpha_2),\\
s_3:(*) \rightarrow &(f_0,f_1,f_2+\frac{\alpha_3 g_1}{f_3},f_3,f_4,g_1,g_2+\frac{\alpha_3}{f_3};\alpha_0,\alpha_1,\alpha_2+\alpha_3,-\alpha_3,\alpha_4),\\
s_4:(*) \rightarrow &(f_0,f_1,f_2+\frac{\alpha_4 g_1}{f_4},f_3,f_4,g_1,g_2+\frac{\alpha_4}{f_4};\alpha_0,\alpha_1,\alpha_2+\alpha_4,\alpha_3,-\alpha_4).
\end{split}
\end{align}
\end{Theorem}

\begin{Theorem}
This system has two first integrals\rm{:\rm}
$$
\frac{d(f_0-f_1)}{dt}=\frac{d(f_3-f_4)}{dt}=0, \quad \frac{d(f_2-g_1g_2)}{dt}=f_2-g_1g_2.
$$
\end{Theorem}
From this, we have
$$
f_0=f_1-1, \quad f_3=f_4-1, \quad f_2-g_1g_2=e^{(t+c)}.
$$
Here we set
$$
t+c=log T, \quad x:=g_1, \ y:=f_1, \ z:=g_2, \ w:=f_4,
$$
then we obtain a 4-parameter family of coupled Painlev\'e III systems in dimension four with affine Weyl group symmetry of type $D_4^{(1)}$ explicitly given by
\begin{equation}
  \left\{
  \begin{aligned}
   \frac{dx}{dT} &=\frac{2x^2y-x^2+(\alpha_0+\alpha_1)x}{T}-1+2w,\\
   \frac{dy}{dT} &=\frac{-2xy^2+2xy-(\alpha_0+\alpha_1)y+\alpha_1}{T},\\
   \frac{dz}{dT} &=\frac{2z^2w-z^2+(\alpha_3+\alpha_4)z}{T}-1+2y,\\
   \frac{dw}{dT} &=\frac{-2zw^2+2zw-(\alpha_3+\alpha_4)w+\alpha_4}{T}\\
   \end{aligned}
  \right. 
\end{equation}
with the Hamiltonian
\begin{align}
\begin{split}
H&=\frac{x^2y^2-x^2y+(\alpha_0+\alpha_1)xy-\alpha_1x}{T}-y\\
&+\frac{z^2w^2-z^2w+(\alpha_3+\alpha_4)zw-\alpha_4z}{T}-w+2yw.
\end{split}
\end{align}

\begin{figure}
\unitlength 0.1in
\begin{picture}(23.71,20.35)(22.18,-23.65)
%
\special{pn 20}%
\special{ar 2462 834 244 244  0.0000000 6.2831853}%
%
\special{pn 20}%
\special{ar 2473 2121 244 244  0.0000000 6.2831853}%
%
\special{pn 20}%
\special{ar 3397 1494 244 244  0.0000000 6.2831853}%
%
\special{pn 20}%
\special{ar 4334 814 244 244  0.0000000 6.2831853}%
%
\special{pn 20}%
\special{ar 4345 2101 244 244  0.0000000 6.2831853}%
%
\special{pn 20}%
\special{pa 2660 988}%
\special{pa 3210 1307}%
\special{fp}%
%
\special{pn 20}%
\special{pa 2638 1934}%
\special{pa 3221 1681}%
\special{fp}%
%
\special{pn 20}%
\special{pa 3573 1318}%
\special{pa 4189 1021}%
\special{fp}%
%
\special{pn 20}%
\special{pa 3551 1681}%
\special{pa 4112 1945}%
\special{fp}%
\put(22.7500,-22.0900){\makebox(0,0)[lb]{$y-1$}}%
\put(23.0800,-9.3300){\makebox(0,0)[lb]{$y$}}%
\put(31.6000,-15.6000){\makebox(0,0)[lb]{$xz+T$}}%
\put(42.0000,-9.2000){\makebox(0,0)[lb]{$w$}}%
\put(41.6000,-21.9000){\makebox(0,0)[lb]{$w-1$}}%
\put(24.2000,-18.6000){\makebox(0,0)[lb]{$0$}}%
\put(24.2000,-5.3000){\makebox(0,0)[lb]{$1$}}%
\put(33.7000,-12.1000){\makebox(0,0)[lb]{$2$}}%
\put(43.2000,-18.4000){\makebox(0,0)[lb]{$3$}}%
\put(43.0000,-5.0000){\makebox(0,0)[lb]{$4$}}%
\end{picture}%
\label{D411}
\caption{The transformations $s_i$ satisfy the relations: $s_i^2=1 \ (i=0,1,2,3,4), \ (s_0s_1)^2=(s_0s_3)^2=(s_0s_4)^2=(s_1s_3)^2=(s_1s_4)^2=(s_3s_4)^2=1, \ (s_0s_2)^3=(s_1s_2)^3=(s_3s_2)^3=(s_4s_2)^3=1.$}
\end{figure}
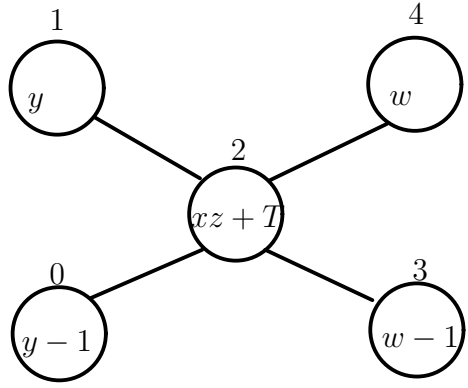

\begin{Theorem}
This system is invariant under the transformations $s_0,\dots,s_4,\pi_1,\\
\pi_2,\pi_3$ defined as follows$:$ with {\it the notation} $(*):=(x,y,z,w,T;\alpha_0,\alpha_1,\alpha_2,\alpha_3,\alpha_4),$
\begin{align}
\begin{split}
s_0:(*) \rightarrow &(x+\frac{\alpha_0}{y-1},y,z,w,T;-\alpha_0,\alpha_1,\alpha_2+\alpha_0,\alpha_3,\alpha_4),\\
s_1:(*) \rightarrow &(x+\frac{\alpha_1}{y},y,z,w,T;\alpha_0,-\alpha_1,\alpha_2+\alpha_1,\alpha_3,\alpha_4),\\
s_2:(*) \rightarrow &(x,y-\frac{\alpha_2 z}{xz+T},z,w-\frac{\alpha_2 x}{xz+T},T;\\
&\alpha_0+\alpha_2,\alpha_1+\alpha_2,-\alpha_2,\alpha_3+\alpha_2,\alpha_4+\alpha_2),\\
s_3:(*) \rightarrow &(x,y,z+\frac{\alpha_3}{w-1},w,T;\alpha_0,\alpha_1,\alpha_2+\alpha_3,-\alpha_3,\alpha_4),\\
s_4:(*) \rightarrow &(x,y,z+\frac{\alpha_4}{w},w,T;\alpha_0,\alpha_1,\alpha_2+\alpha_4,\alpha_3,-\alpha_4),\\
\pi_1:(*) \rightarrow &(-x,1-y,z,w,-T;\alpha_1,\alpha_0,\alpha_2,\alpha_3,\alpha_4),\\
\pi_2:(*) \rightarrow &(x,y,-z,1-w,-T;\alpha_0,\alpha_1,\alpha_2,\alpha_4,\alpha_3),\\
\pi_3:(*) \rightarrow &(z,w,x,y,T;\alpha_3,\alpha_4,\alpha_2,\alpha_0,\alpha_1).
\end{split}
\end{align}
\end{Theorem}

\end{document}